\newif\ifIAI
\IAIfalse

\ifIAI
\documentclass[numbers,webpdf,imaiai]{ima-authoring-template}
\theoremstyle{thmstyletwo}%
\else
\documentclass[letterpaper]{article}
\usepackage[margin=1in]{geometry}
\usepackage[numbers]{natbib}
\usepackage{amssymb}
\usepackage{amsthm}
\usepackage{mathrsfs}
\usepackage{authblk}

\usepackage{hyperref}
\hypersetup{
	colorlinks,
	linkcolor=blue,
	citecolor=blue,
	urlcolor=blue,
}
\fi

\usepackage{import}
\usepackage{mathtools}

\usepackage{xparse}
\usepackage{import}
\usepackage[normalem]{ulem}

\makeatletter
\@ifpackageloaded{unicode-math}{%
	\newcommand{\mymathbold}{\symbf}%
}{%
	\usepackage{bm}%
	\newcommand{\mymathbold}{\bm}%
}
\makeatother

\newcommand{\scrbar}[1]{\overline{\mathcal{#1}}}
\newcommand{\scrhat}[1]{\widehat{\mathcal{#1}}}
\newcommand{\scrtl}[1]{\widetilde{\mathcal{#1}}}

\ExplSyntaxOn

\cs_new_protected:Nn \bamboo_define:nnnnnN
{
	\cs_new_protected:cpx { #3 #1 #4 } { \exp_not:N #5{#6{#2}} }
}

\int_step_inline:nnn { `A } { `Z }
{
	\bamboo_define:nnnnnN
	{ \char_generate:nn { #1 } { 11 } }
	{ \char_generate:nn { #1 } { 11 } }
	{ bl }
	{ }
	{ \mymathbold }
	\use:n
	
	\bamboo_define:nnnnnN
	{ \char_generate:nn { #1 } { 11 } }
	{ \char_generate:nn { #1 } { 11 } }
	{ scr }
	{ }
	{ \mathcal }
	\use:n
	
	\bamboo_define:nnnnnN
	{ \char_generate:nn { #1 } { 11 } }
	{ \char_generate:nn { #1 } { 11 } }
	{ }
	{ hat }
	{ \widehat }
	\use:n
	
	\bamboo_define:nnnnnN
	{ \char_generate:nn { #1 } { 11 } }
	{ \char_generate:nn { #1 } { 11 } }
	{ scr }
	{ hat }
	{ \scrhat }
	\use:n
	
	\bamboo_define:nnnnnN
	{ \char_generate:nn { #1 } { 11 } }
	{ \char_generate:nn { #1 } { 11 } }
	{ }
	{ br }
	{ \overline }
	\use:n
	
	\bamboo_define:nnnnnN
	{ \char_generate:nn { #1 } { 11 } }
	{ \char_generate:nn { #1 } { 11 } }
	{ }
	{ tl }
	{ \widetilde }
	\use:n
	
	\bamboo_define:nnnnnN
	{ \char_generate:nn { #1 } { 11 } }
	{ \char_generate:nn { #1 } { 11 } }
	{ scr }
	{ br }
	{ \scrbar }
	\use:n
	
	\bamboo_define:nnnnnN
	{ \char_generate:nn { #1 } { 11 } }
	{ \char_generate:nn { #1 } { 11 } }
	{ scr }
	{ tl }
	{ \scrtl }
	\use:n
	
	\bamboo_define:nnnnnN
	{ \char_generate:nn { #1 } { 11 } }
	{ \char_generate:nn { #1 } { 11 } }
	{ }
	{ dt }
	{ \dot}
	\use:n
}
\int_step_inline:nnn { `a } { `z }
{
	\bamboo_define:nnnnnN
	{ \char_generate:nn { #1 } { 11 } }
	{ \char_generate:nn { #1 } { 11 } }
	{ bl }
	{ }
	{ \mymathbold }
	\use:n
	
	\bamboo_define:nnnnnN
	{ \char_generate:nn { #1 } { 11 } }
	{ \char_generate:nn { #1 } { 11 } }
	{ }
	{ hat }
	{ \hat }
	\use:n
	
	\bamboo_define:nnnnnN
	{ \char_generate:nn { #1 } { 11 } }
	{ \char_generate:nn { #1 } { 11 } }
	{ }
	{ br }
	{ \bar }
	\use:n
	
	\bamboo_define:nnnnnN
	{ \char_generate:nn { #1 } { 11 } }
	{ \char_generate:nn { #1 } { 11 } }
	{ }
	{ tl }
	{ \tilde }
	\use:n                            
	
	\bamboo_define:nnnnnN
	{ \char_generate:nn { #1 } { 11 } }
	{ \char_generate:nn { #1 } { 11 } }
	{ }
	{ dt }
	{ \dot}
	\use:n
}
\clist_map_inline:nn
{
	Gamma,Delta,Theta,Lambda,Xi,Pi,Sigma,Phi,Psi,Omega
}
{
	\bamboo_define:nnnnnN
	{ #1 }
	{ #1 }
	{ bl }
	{ }
	{ \mymathbold }
	\use:c
	
	\bamboo_define:nnnnnN
	{ #1 }
	{ #1 }
	{ op }
	{ }
	{ \op }
	\use:c
	
	\bamboo_define:nnnnnN
	{ #1 }
	{ #1 }
	{ scr }
	{ }
	{ \mathcal }
	\use:c
	
	\bamboo_define:nnnnnN
	{ #1 }
	{ #1 }
	{ }
	{ hat }
	{ \widehat }
	\use:c
	
	\bamboo_define:nnnnnN
	{ #1 }
	{ #1 }
	{ }
	{ br }
	{ \overline }
	\use:c
	
	\bamboo_define:nnnnnN
	{ #1 }
	{ #1 }
	{ }
	{ tl }
	{ \widetilde }
	\use:c
	
	\bamboo_define:nnnnnN
	{ #1 }
	{ #1 }
	{ }
	{ dt }
	{ \dot}
	\use:c
	
}
\clist_map_inline:nn
{
	alpha,beta,gamma,delta,epsilon,zeta,eta,theta,iota,kappa,
	lambda,mu,nu,xi,pi,rho,sigma,tau,phi,chi,psi,omega
}
{
	\bamboo_define:nnnnnN
	{ #1 }
	{ #1 }
	{ bl }
	{ }
	{ \mymathbold }
	\use:c
	
	\bamboo_define:nnnnnN
	{ #1 }
	{ #1 }
	{ }
	{ hat }
	{ \hat }
	\use:c
	
	\bamboo_define:nnnnnN
	{ #1 }
	{ #1 }
	{ }
	{ br }
	{ \bar }
	\use:c
	
	\bamboo_define:nnnnnN
	{ #1 }
	{ #1 }
	{ }
	{ tl }
	{ \tilde }
	\use:c
	
	\bamboo_define:nnnnnN
	{ #1 }
	{ #1 }
	{ }
	{ dt }
	{ \dot}
	\use:c
}

\ExplSyntaxOff

\DeclareMathOperator{\E}{\mathbf{E}}
\renewcommand{\P}{\operatorname{\mathbf{P}}}

\newcommand{\trace}{\operatorname{trace}}
\newcommand{\tr}{\operatorname{tr}}

\newcommand{\diag}{\operatorname{diag}}

\DeclarePairedDelimiter{\norm}{\lVert}{\rVert}

\DeclarePairedDelimiter{\abs}{\lvert}{\rvert}

\DeclarePairedDelimiter{\parens}{(}{)}
\DeclarePairedDelimiter{\brackets}{[}{]}
\DeclarePairedDelimiterX{\ip}[2]{\langle}{\rangle}{#1,#2}

\DeclarePairedDelimiterXPP{\normsub}[2]{}{\lVert}{\rVert}{_{#2}}{#1}
\DeclarePairedDelimiterXPP{\ipsub}[3]{}{\langle}{\rangle}{_{#3}}{#1,#2}

\DeclarePairedDelimiterXPP{\ipHS}[2]{}{\langle}{\rangle}{_{\mathrm{HS}}}{#1, #2}
\DeclarePairedDelimiterXPP{\normHS}[1]{}{\lVert}{\rVert}{_{\mathrm{HS}}}{#1}

\DeclarePairedDelimiterXPP{\ipF}[2]{}{\langle}{\rangle}{_{\mathrm{F}}}{#1, #2}
\DeclarePairedDelimiterXPP{\normF}[1]{}{\lVert}{\rVert}{_{\mathrm{F}}}{#1}

\DeclarePairedDelimiterXPP{\dkl}[2]{\operatorname{D_{KL}}}{(}{)}{}{#1 \: \delimsize\Vert \: #2}

\DeclarePairedDelimiterXPP{\restr}[2]{}{{}}{\vert}{_{#2}}{#1}

\newcommand{\R}{\mathbf{R}}

\newcommand{\Z}{\mathbf{Z}}

\newcommand{\spn}{\operatorname{span}}

\newcommand{\binomdist}{\operatorname{Binom}}

\newcommand{\todo}[1]{\textbf{[TODO: #1]}}

\usepackage{bookmark}
\usepackage{color}
\usepackage[capitalize]{cleveref}
\usepackage[left]{lineno}

\newtheorem{theorem}{Theorem}
\newtheorem{corollary}{Corollary}
\newtheorem{lemma}{Lemma}

\DeclarePairedDelimiterXPP{\opnorm}[1]{}{\lVert}{\rVert}{_{\mathrm{op}}}{#1}
\DeclarePairedDelimiterXPP{\nucnorm}[1]{}{\lVert}{\rVert}{_{*}}{#1}
\DeclarePairedDelimiterXPP{\normFtr}[1]{}{\lVert}{\rVert}{_{\mathrm{F,tr}}}{#1}
\DeclarePairedDelimiterXPP{\ipFtr}[2]{}{\langle}{\rangle}{_{\mathrm{F,tr}}}{#1, #2}
\newcommand\ddiag{\operatorname{ddiag}}

\newcommand{\Ogrp}{\mathrm{O}}

\newcommand{\st}{\text{ s.t.\ }}

\newcommand{\onevec}{\mathbf{1}}
\newcommand{\ER}{Erd\H{o}s--R\'enyi}
\newcommand{\ERmath}{\mathcal{G}}

\newcommand{\Pone}{P_\onevec}
\newcommand{\Poneperp}{P_{\onevec}^\perp}

\ifIAI

\linenumbers
\renewcommand{\todo}[1]{}
\fi

\newcommand{\longtitletext}{Nonconvex landscapes for $\mathbf{Z}_2$ synchronization and graph clustering \\are benign near exact recovery thresholds}
\newcommand{\shorttitletext}{Nonconvex landscapes for $\Z_2$ synchronization and graph clustering}
\newcommand{\abstracttext}{%
We study the optimization landscape of a smooth nonconvex program arising from synchronization over the two-element group $\mathbf{Z}_2$, that is, recovering $z_1, \dots, z_n \in \{\pm 1\}$ from (noisy) relative measurements $R_{ij} \approx z_i z_j$. Starting from a max-cut--like combinatorial problem, for integer parameter $r \geq 2$, the nonconvex problem we study can be viewed both as a rank-$r$ Burer--Monteiro factorization of the standard max-cut semidefinite relaxation and as a relaxation of $\{ \pm 1 \}$ to the unit sphere in $\mathbf{R}^r$. First, we present deterministic, non-asymptotic conditions on the measurement graph and noise under which every second-order critical point of the nonconvex problem yields exact recovery of the ground truth. Then, via probabilistic analysis, we obtain asymptotic guarantees for three benchmark problems: (1) synchronization with a complete graph and Gaussian noise, (2) synchronization with an Erdős--Rényi random graph and Bernoulli noise, and (3) graph clustering under the binary symmetric stochastic block model. In each case, we have, asymptotically as the problem size goes to infinity, a benign nonconvex landscape near a previously-established optimal threshold for exact recovery; we can approach this threshold to arbitrary precision with large enough (but finite) rank parameter $r$. In addition, our results are robust to monotone adversaries.%
}

\begin{document}
\ifIAI
\else
\title{\longtitletext}
\author[1]{Andrew D.\ McRae}
\author[2]{Pedro Abdalla}
\author[2]{Afonso S.\ Bandeira}
\author[1]{Nicolas Boumal}
\affil[1]{Institute of Mathematics, EPFL\authorcr \texttt{\{andrew.mcrae,nicolas.boumal\}@epfl.ch}}
\affil[2]{Department of Mathematics, ETH Zurich\authorcr \texttt{pedro.abdallateixeira@ifor.math.ethz.ch}\authorcr\texttt{bandeira@math.ethz.ch}}

\maketitle

\abstract{\abstracttext}

\fi

\section{Introduction and result highlights}
\label{sec:intro}
In a \emph{synchronization} problem, one wants to estimate $n$ group elements $g_1, \ldots, g_n$ from pairwise relative measurements $R_{ij} \approx g_i g_j^{-1}$.
Such problems are important in robotics \cite{Rosen2019}, computer vision \cite{Hartley2013}, signal processing \cite{Cucuringu2012,Iwen2020}, computer science \cite{Pachauri2013,Abbe2018}, and many other areas.
In this paper, we study synchronization over the two-element group $\Z_2$ formulated as the real orthogonal group $\Ogrp(1) = \{\pm 1\}$ under multiplication.

The $\Z_2$ synchronization approach we consider takes the form of a combinatorial optimization problem (arising from a maximum likelihood and/or least squares formulation of the synchronization problem---see, e.g., \cite{Singer2011,Bandeira2017}):
\begin{equation}
	\label{eq:obj_orig}
	\max_{x\in \{\pm 1\}^n}~\sum_{i,j=1}^n C_{ij} x_i x_j
	= \max_{x\in \{\pm 1\}^n}~\ip{C}{x x^\top},
\end{equation}
where $C \in \R^{n\times n}$ is a cost matrix to be specified based on the exact application,
and $\ip{A}{B} = \trace(A^\top B)$ denotes the entrywise Euclidean (Frobenius) inner product between real matrices $A$ and $B$ of the same size.
As the NP-hard max-cut problem is an instance of \eqref{eq:obj_orig},
we expect the problem to be computationally intractable in general
(though not necessarily for particular instances; indeed, in the synchronization context, our results will yield a tractable algorithm for computing the exact solution).
A classical (approximate) solution approach is the semidefinite program (SDP) relaxation of \citet{Goemans1995}:
\begin{equation}
	\label{eq:obj_sdp}
	\max_{X \succeq 0}~\ip{C}{X} \st \ddiag(X) = I_n,
\end{equation}
where $\ddiag \colon \R^{n \times n} \to \R^{n \times n}$ extracts the diagonal part of a matrix, and $I_n$ is the $n\times n$ identity matrix.
For synchronization problems, this has proved to be a powerful and, in many cases, statistically optimal approach \cite{Hajek2016a,Bandeira2017,Gao2022}.
However, as we have approximately squared the number of variables compared to the original problem \eqref{eq:obj_orig}, it is still computationally expensive to solve \eqref{eq:obj_sdp} directly for large $n$.

To lighten the computational burden of the SDP approach, we consider its Burer--Monteiro factorization \cite{Burer2003},
which corresponds to the following family of smooth but nonconvex ``partial'' relaxations of \eqref{eq:obj_orig}: for any integer $r \geq 2$,
\begin{equation}
	\label{eq:obj_bm}
	\max_{Y \in \R^{n\times r}}~\ip{C}{YY^\top} \st \ddiag(YY^\top) = I_n.
\end{equation}
Although \eqref{eq:obj_bm} is nonconvex, in many cases it has a \emph{benign landscape} in the sense that every second-order critical point $Y$ (where the gradient is zero and the Hessian is negative semidefinite in a Riemannian geometric sense) is in fact globally optimal and also solves the SDP relaxation (i.e., $X = Y Y^\top$ solves \eqref{eq:obj_sdp}).
This is generically true when $r \approx \sqrt{n}$ \cite{Boumal2019}.
However, for the specific problem instances arising in synchronization,
there is much empirical and theoretical work showing that \eqref{eq:obj_bm} has a benign landscape even for small $r$ not increasing with $n$ \cite{Bandeira2016a,Ling2023b,Ling2023a,McRae2024};
in this case, we can solve \eqref{eq:obj_sdp} via a local algorithm on the problem \eqref{eq:obj_bm}, which has only $rn$ variables.
See, for example, \citet{Ling2023a,McRae2024} for an introduction to this approach for general orthogonal group synchronization problems.

Our main deterministic result (\Cref{thm:z2_determ}) provides sufficient conditions on the $\Z_2$ synchronization problem (in terms of the noise level and measurement graph connectivity) for which not only does the corresponding problem \eqref{eq:obj_bm} have a benign landscape for small $r$,
but also the solutions \emph{exactly} recover the underlying ground-truth signal.\footnote{More precisely, we recover the ground truth up to a global sign to which the problem is invariant.}
Our conditions are sharper than previous such results.

We then combine this deterministic result with probabilistic analyses to obtain asymptotic exact recovery results for several common random models in $\Z_2$ synchronization and graph clustering.
In each case, our results match established optimal thresholds for exact recovery.
This was previously only shown for the direct SDP approach or other specialized algorithms;
previous analyses of the nonconvex approach fell short of the optimal thresholds.
Furthermore, adapting the techniques of \citet{Ling2023b}, we show that our results largely inherit one of the key advantages of the SDP approach, which is its robustness to \emph{monotone adversaries} (perturbations which seem ``helpful'' but can disrupt algorithms by changing the problem's random structure).

Before describing our results in full generality, we now present in brief these key statistical exact recovery guarantees.

\subsection{$\Z_2$ synchronization with Gaussian noise.}
\label{sec:intro_gauss}
As a simple example and application of our results,
consider the problem of recovering a vector $z\in \{\pm 1\}^n$ from the noisy relative measurements $R_{ij} = z_i z_j + \sigma W_{ij}$ (symmetric in $i,j$),
where $\sigma \geq 0$ controls the noise strength,
and $W_{ij}$ are independent and identically distributed (i.i.d.)\ standard normal random variables (modulo symmetry: $W_{ij} = W_{ji}$). Up to rescaling (and the choice of diagonal elements, which have no effect on our problem), $W$ is often called a Gaussian Wigner matrix or a draw from the Gaussian Orthogonal Ensemble.
\citet{Bandeira2017} introduced this particular model as a simpler version of the well-known angular synchronization problem. \citet{Javanmard2016} later used it as an alternative to the Bernoulli-noise problem we consider in \Cref{sec:intro_er}.

We have the following asymptotic recovery result for this model (proved in \Cref{sec:z2_proof_gauss}):
\begin{corollary}
	\label{cor:z2_gauss}
	For $n \geq 2$, let
	\begin{equation}
		\label{eq:z2_gauss_model}
		C = z z^\top + \sigma W \in \R^{n \times n},
	\end{equation}where $W$ is a symmetric $n \times n$ random matrix with i.i.d.\ (up to symmetry) standard normal off-diagonal entries.
	Fix a tolerance $\epsilon>0$ and an integer $r_0 \geq 3$ independent of $n$, and suppose
	\begin{equation*}
		\sigma \leq \frac{r_0-3}{r_0-1}\sqrt{\frac{n}{(2+\epsilon)\log n}}.
	\end{equation*}
	Then, with probability $\to 1$ as $n \to \infty$, for all $r \geq r_0$, every second-order critical point $Y$ of \eqref{eq:obj_bm} with cost matrix $C$ satisfies $Y = z u^\top$ for some unit-norm $u \in \R^r$.
\end{corollary}
From this, we can recover $z$ up to global sign as the top left singular vector of $Y$.

\citet{Bandeira2017,Bandeira2018} showed that the SDP relaxation achieves exact recovery when $\sigma \leq \sqrt{\frac{n}{(2+\epsilon)\log n}}$
and that this is the best possible threshold in terms of the noise variance~$\sigma$.
Our result shows that nonconvex optimization of \eqref{eq:obj_bm} achieves the same threshold for large enough $r$ (depending only on $\epsilon$);
thus we can approach the optimal recovery threshold by nonconvex optimization over $\sim C(\epsilon) n$ variables rather than $\sim n^2$ as with the semidefinite relaxation \eqref{eq:obj_sdp}.

Previously, \citet{Bandeira2016a} showed that the nonconvex problem \eqref{eq:obj_bm} with $r = 2$
can achieve exact recovery but with suboptimal scaling of $\sigma$ in the problem size $n$.
More recently, \citet{Ling2023b} proved this result with the optimal scaling, showing that for large enough $r$, $\sigma \leq \frac{1}{4 + \epsilon} \sqrt{\frac{n}{\log n}}$ suffices.
\Cref{cor:z2_gauss} improves the constant to reach the exact asymptotic threshold.

\subsection{$\Z_2$ synchronization with partial measurements and Bernoulli noise}
\label{sec:intro_er}
The model \eqref{eq:z2_gauss_model} assumes that we have a measurement for \emph{every} pair $(i,j)$, which is too restrictive for many applications.
More generally, we may only have measurements for certain pairs $(i,j)$.
We model this by a measurement graph $G=(V,E)$, where $V = \{1, \dots, n\}$, and we only have access to measurements $R_{ij} \approx z_i z_j$ for pairs $(i,j) \in E$.

For now, we consider a simple random graph model (we study more general deterministic graphs in \Cref{sec:z2_deterministic}).
Suppose the measurement graph $G = (V, E)$ is sampled from the \ER{} random graph distribution $\ERmath(n, p)$,
in which each possible edge $(i,j)$ appears with fixed probability $p$ (independent of the others).
Clearly, if the graph $G$ is disconnected, then global synchronization is impossible. The \ER{} random graph has a connectivity phase transition at $p = \frac{\log n}{n}$.
Thus we will consider $p=\frac{a \log n}{n}$ for some constant parameter $a>1$ because only in this regime is $G$ connected with high probability for large $n$.
We consider more general scaling of $p$ with $n$ in \Cref{sec:z2_bern}.

Furthermore, although the Gaussian noise in \eqref{eq:z2_gauss_model} is convenient for analysis,
it may not be the most natural noise model. Indeed, we have the prior knowledge that the entries of $z$ (and hence the ideal observations of them) are $\pm 1$.
This motivates another \emph{multiplicative} noise model in which the signs of our observations are flipped with some probability. We assume that the noise acts independently and identically (i.e., each observed sign is flipped with a constant probability) on each pair $(i,j)$:
this is the Bernoulli noise model.

Specifically, given a sample $G = (V, E)$ from $\ERmath(n,p)$,
we want to recover a vector $z\in \{\pm 1\}^n$ from the measurements
\begin{equation}
	\label{eq:z2_bernoulli_model}
	R_{ij} = 
	\begin{cases}
		z_i z_j & \text{with probability } \frac{1 + \delta}{2} \\
		- z_i z_j & \text{with probability } \frac{1 - \delta}{2}
	\end{cases}
	\quad \text{for $(i, j) \in E$,}
\end{equation}
where $\delta \in [0, 1]$ determines how likely we are to see the correct sign.
Larger $\delta$ implies better measurements;
if $\delta = 0$, the measurements are independent of the signal $z$ (and therefore useless),
whereas $\delta = 1$ implies our measurements are entirely uncorrupted.

Note that we can also interpret this problem as clustering graph vertices based on \emph{signed} edges:
if $z_i$ is the ``true'' cluster label of node $i$, an edge between nodes $i$ and $j$ is likely to be positive (resp.\ negative) if $i$ and $j$ have the same (resp.\ different) label.
This general problem arose in a machine learning context as ``correlation clustering'' \cite{Bansal2004}.
The ``$\Z_2$ synchronization'' interpretation was introduced by \citet{Cucuringu2012,Cucuringu2012a}.
\citet{Abbe2014} proposed the specific probabilistic model we present here.

Our next result establishes when we can obtain exact recovery of $z$ with nonconvex optimization
depending on the graph connectivity and signal strength parameters:
\begin{corollary}
	\label{cor:z2_bern}
	Consider the model \eqref{eq:z2_bernoulli_model} with $\delta \in [0, 1]$, where the graph $G$ is sampled from $\mathcal{G}(n,p)$ for $p=\frac{a \log n}{n}$ with fixed $a > 1$.
	Suppose, for some $\epsilon>0$,
	\[
		a (1 - \sqrt{1 - \delta^2}) \geq 1 + \epsilon.
	\]
	Then there exists a constant $r_0$ (depending only on $a$ and $\epsilon$) such that, with probability $\to 1$ as $n \to \infty$, the following holds:
	if we choose the cost matrix $C$ according to
	\begin{equation}
		\label{eq:bern_er_C}
		C_{ij} = \begin{cases}
			R_{ij} & \text{if } (i,j) \in E \\
			0 & \text{otherwise},
		\end{cases}
	\end{equation}
	then, for all $r \geq r_0$, every second-order critical point $Y$ of \eqref{eq:obj_bm} satisfies $Y = z u^\top$ for some unit-norm $u \in \R^r$.
\end{corollary}
This threshold is the best possible for exact recovery \cite{Hajek2016a,Jog2015}.
In particular, \citet{Hajek2016a} showed that this threshold is achieved by the semidefinite relaxation \eqref{eq:obj_sdp};
once again, we show that we can reach the same threshold via nonconvex optimization of \eqref{eq:obj_bm}.
\Cref{cor:z2_bern} is a simplified version of the more general result \Cref{thm:z2_er_bern} in \Cref{sec:z2_bern}; see the discussion after that theorem for further context and related work.

\subsection{Graph clustering with the binary stochastic block model}
\label{sec:intro_sbm}
We now shift focus to a different measurement model.
We are again given a graph and we want to cluster its nodes.
In \Cref{sec:intro_er}, we had signed edges indicating whether nodes are in the same cluster,
but now we have no signs and must work only with edge density. We want to assign the nodes to clusters such that the density of edges within clusters is larger than the density between clusters.
This is a classical and well-studied problem and is also called ``community detection'' or the ``planted partition'' (or bisection) problem.
See the surveys of \citet{Fortunato2016,Moore2017,Abbe2018} and the citations below for further introduction.
We focus on one specific setup in this area.

The \emph{stochastic block model} (SBM) for a graph is the following:
given ``ground truth'' clusters,
we assume that nodes in the same cluster are connected with an edge with some fixed probability $p$,
and that nodes in different clusters are connected with fixed probability $q < p$
(the existence of each edge being independent of the others).

In particular, we consider the binary and symmetric case, in which there are two equally-sized clusters.
Let $n$ be even,
and consider the following random graph $G$ on vertices $1, \dots, n$.
Given a vector $z \in \{ \pm 1 \}^n$ of ``true'' cluster labels (with $\ip{z}{\onevec} = 0$, as the cluster have equal size; $\onevec$ denotes the all-ones vector),
the entries $A_{ij}$ of the adjacency matrix of $G$ are independent (modulo symmetry) Bernoulli random variables with, for $i \neq j$ (we assume $A$ has zero diagonal),
\begin{equation}
	\label{eq:sbm_model}
	\P(A_{ij} = 1) = \begin{cases}
		p & \text{if } z_i = z_j, \text{ and}\\
		q & \text{if } z_i \neq z_j.
	\end{cases}
\end{equation}
The connection to synchronization becomes clear from the fact that
\begin{equation}
	\label{eq:EA}
	\E A = \frac{p - q}{2} z z^\top + \frac{p + q}{2} \onevec \onevec^\top - p I_n.
\end{equation}
The last term reflects the fact that the diagonal of $A$ is zero.

However, treating this exactly like $\Z_2$ synchronization by taking cost matrix $C = A$ in \eqref{eq:obj_orig} would not work, because the all-ones vector is a trivial solution.
This is due to the dominance of the all-ones term in \eqref{eq:EA}.
We can deal with this by (approximately) subtracting off the mean of the entries of $A_{ij}$:
we set
\begin{equation}
	\label{eq:sbm_C}
	C_{ij} \coloneqq A_{ij} - \frac{1}{n^2} \sum_{i,j=1}^n A_{ij}, \quad \text{that is,} \quad
	C \coloneqq A - \frac{1}{n^2} \ip{A}{\onevec \onevec^\top} \onevec \onevec^\top.
\end{equation}
If the true average edge probability is known, we can more simply set
\begin{equation}
	\label{eq:sbm_C_pqknown}
	C \coloneqq A - \frac{p + q}{2} \onevec \onevec^\top.
\end{equation}
Because this graph resembles the \ER{} random graph considered in the previous section,
it is natural once again to set the edge probabilities proportional to $\frac{\log n}{n}$.
We then have the following result:
\begin{corollary}
	\label{cor:SBM_BM}
	Consider the model \eqref{eq:sbm_model} with parameters $p = \frac{a \log n}{n}$ and $q = \frac{b \log n}{n}$
	for some $a, b \geq 0$ satisfying $\sqrt{a} - \sqrt{b} > \sqrt{2}$.
	There exists a constant $r_0$ (depending only on $a$ and $b$) such that, with probability $\to 1$ as $n \to \infty$, the following holds:
	if we choose the cost matrix $C$ according to \eqref{eq:sbm_C} or \eqref{eq:sbm_C_pqknown},
	then, for all $r \geq r_0$, every second-order critical point $Y$ of \eqref{eq:obj_bm} satisfies $Y = z u^\top$ for some unit-norm $u \in \R^r$.
\end{corollary}
The condition on $a$ and $b$ is optimal for exact recovery \cite{Mossel2015,Abbe2016}.
Once again, we have shown that optimal results for the SDP relaxation \eqref{eq:obj_sdp} (found in \cite{Hajek2016,Bandeira2018}) can be achieved with nonconvex optimization in \eqref{eq:obj_bm}.
This is a simplified version of the more general result \Cref{thm:sbm_asymp} in \Cref{sec:sbm_main};
see the discussion after that theorem for details and comparison to previous work.

\subsection{Paper outline and notation}
The remainder of this paper is organized as follows:
\begin{itemize}
	\item \Cref{sec:z2_deterministic} presents a general deterministic result (\Cref{thm:z2_determ}) that describes when, in the context of $\Z_2$ synchronization, the nonconvex problem \eqref{eq:obj_bm} has a benign landscape such that all second-order critical points yield exact recovery of the ground truth.
	
	\item \Cref{sec:z2_asymp} uses this deterministic result together with probabilistic analyses to study the random $\Z_2$ synchronization problems presented in \Cref{sec:intro_gauss,sec:intro_er}.
	In addition to \Cref{cor:z2_gauss} for Gaussian noise, we obtain \Cref{thm:z2_er_bern}, which is a more general version of \Cref{cor:z2_bern} for \ER{} random graphs and Bernoulli noise.
	
	\item \Cref{sec:sbm_main} adapts these arguments to the graph clustering problem presented in \Cref{sec:intro_sbm};
	we obtain the high-probability asymptotic result \Cref{thm:sbm_asymp} from which \Cref{cor:SBM_BM} is derived.
	Along the way, we obtain a deterministic condition (\Cref{cor:sbm_determ}) that may be of independent interest.
	
	\item \Cref{sec:monotone} discusses the robustness (summarized in \Cref{thm:monotone}) of our results to monotone adversaries.
\end{itemize}

In general, we define notation as needed when we use it.
However, we include here for completeness some basic definitions that are common in the literature and will appear throughout the paper.
	$\onevec$ denotes the all-ones vector in $\R^n$.
	$\norm{v}$ denotes the Euclidean ($\ell_2$) norm of a vector $v$.
	$\opnorm{M}$, $\normF{M}$, and $\nucnorm{M}$ respectively denote the operator (in $\ell_2$) norm, Frobenius (elementwise $\ell_2$) norm, and nuclear norm (sum of singular values) of a matrix $M$.
	$\ip{\cdot}{\cdot}$ is the inner product notation both for the standard Euclidean inner product of vectors and the elementwise (Frobenius) inner product of matrices (the distinction will be clear from context).
	$\diag(v)$, for a vector $v \in \R^n$, denotes the $n \times n$ diagonal matrix with the entries of $v$ on the diagonal.
	$\ddiag\colon \R^{n \times n} \to \R^{n \times n}$ extracts the diagonal of a matrix and sets the remaining elements to zero.
	$a \lesssim b$ means there is an unspecified but universal (i.e., not depending on any problem parameters) constant $c > 0$ such that $a \leq c b$ (and $a \gtrsim b$ means $b \lesssim a$).

\section{Deterministic $\Z_2$ synchronization analysis}
\label{sec:z2_deterministic}
In this section, we state and prove our main deterministic result for $\Z_2$ synchronization.
We start by stating the problem in greater generality than before.
Let $A \in \R^{n \times n}$ be the symmetric adjacency matrix of an undirected but potentially weighted graph $G$ on $n$ vertices.
Assume the cost matrix $C$ in \eqref{eq:obj_orig} has the form
\begin{equation}
	\label{eq:z2_model}
	C_{ij} = A_{ij} z_i z_j + \Delta_{ij},\quad \text{that is,} \quad C = \diag(z) A \diag(z) + \Delta,
\end{equation}
where $z=(z_1,\ldots,z_n) \in \{\pm 1\}^n$ is the vector we want to recover, and $\Delta$ is a matrix representing noise.
Notice that the models \eqref{eq:z2_gauss_model} and \eqref{eq:bern_er_C} are both instances of \eqref{eq:z2_model}.

In general, we will assume, without loss of generality, that $C$, $A$, and $\Delta$ have zero diagonal,
as this has no effect on the landscape of \eqref{eq:obj_bm} (only changing the objective value by a constant).

We consider the following question: which conditions on the graph $G$ (i.e, its adjacency matrix $A$) and the noise matrix $\Delta$ ensure that we can recover $z$ via \eqref{eq:obj_bm}?
\begin{itemize}
\item \textbf{Connectivity:} Clearly, if the graph $G$ is disconnected, synchronization is impossible.
Thus we expect some measure of connectivity to be important.
In particular, recall that the Laplacian matrix of $G$ is $L \coloneqq \diag(A \onevec) - A = D-A$,
where $D$ is the $n \times n$ diagonal matrix whose diagonal entries are the vertex degrees $d_1, \dots, d_n$:
\[
	D_{ii} = d_i \coloneqq \sum_{\substack{j = 1\\j \neq i}}^n A_{ij}.
\]
It is well known that $L$ is positive semidefinite with eigenvalues $0=\lambda_1\le \lambda_2 \le \ldots \le \lambda_n$.
Furthermore, $G$ is connected if and only $\lambda_2>0$.
The quantity $\lambda_2$ is called the algebraic connectivity of $G$ and plays a key role in our analysis.

\item \textbf{Local Stability:} Notice that to ensure we recover $z$ exactly even with the original combinatorial problem \eqref{eq:obj_orig}, we must rule out the possibility that flipping one entry of $z$ could increase (or fail to increase) the objective function.
In particular, we must have, for all $i$,
\begin{equation*}
	\ip{C}{z z^\top} > \ip{C}{z_{-i} z_{-i}^\top},
\end{equation*}
where $z_{-i}:=(z_1,\ldots,z_{i-1}, -z_i, z_{i+1}, \ldots,z_n)$. Expanding this via \eqref{eq:z2_model} reveals that this condition is equivalent to
\begin{equation}
	\label{eq:z2_stability_cond}
	d_i > -z_i\sum_{\substack{j = 1\\j \neq i}}^n \Delta_{ij} z_j \eqqcolon\rho_i^{\Delta}.
\end{equation}
Therefore, it is natural to expect that a control on the quantity
\begin{equation}
	\label{eq:rhomax}
	\rho^{\Delta}\coloneqq \max_i~\rho_i^{\Delta}
\end{equation}
will be necessary.

\item\textbf{Noise spectral norm:}
Beyond the local stability condition above, it is also useful to control the noise $\Delta$ via its matrix operator norm $\opnorm{\Delta}$.
\end{itemize}

For the Burer--Monteiro approach \eqref{eq:obj_bm}, our main result shows that a suitable condition on those three quantities (depending also on the factorization rank $r$) leads to exact recovery of $z$:
\begin{theorem}
	\label{thm:z2_determ}
	Let $G$ be a connected graph on $n$ vertices with adjacency matrix $A$ and algebraic connectivity $\lambda_2$. Let $C$ be the cost matrix of the $\Z_2$ synchronization measurement model \eqref{eq:z2_model} on $G$. Let $r \geq 4$ be an integer, and suppose
	\begin{equation}
		\label{eq:z2_cond}
		\rho^\Delta + \frac{r+11}{r-3} \opnorm{\Delta} \leq \frac{r-3}{r-1} \lambda_2.
	\end{equation}
	Then every second-order critical point $Y$ of \eqref{eq:obj_bm} satisfies $Y = z u^\top$ for some unit-norm $u \in \R^r$.
	The same result holds if $r = 3$ and $\Delta = 0$.
\end{theorem}
We defer the proof to the end of this section. When $\Delta = 0$, we recover known results \cite{Markdahl2018,McRae2024}.
The result and its proof build on and resemble those of \citet{McRae2024,Ling2023b}.
Compared to \citet{Ling2023b}, our result only depends on $\lambda_2$ rather than a condition number and is thus tighter and more general.
Compared to \citet{McRae2024}, our analysis is specialized to the $\Z_2$ case in order to achieve tighter results in terms of the noise $\Delta$.

More broadly, our work builds on recent results in nonconvex optimization for orthogonal group synchronization (often known, for particular problems, as phase/angular synchronization or synchronization of rotations).
See, for example, \citet{Ling2023a,McRae2024} for more on this topic.
However, most existing results in this area do not directly yield exact recovery guarantees in the presence of noise (in fact, it is not, in general, possible to recover the ground truth exactly with noise because every orthogonal group other than $\Z_2$ is continuous-valued).
Exceptions to this are the results of \citet{Bandeira2016a,Ling2023b} as discussed previously; these papers specifically consider the $\Z_2$ case.
Indirectly, we could combine the more general benign landscape results with other previous work showing that the ground truth is the global optimum of the SDP relaxation (e.g., \cite{Hajek2016,Hajek2016a,Bandeira2018}), but the assumptions required on the problem parameters would be suboptimal.

\Cref{thm:z2_determ} provides \emph{sufficient} conditions for exact recovery, but we make no claim that the conditions are \emph{necessary} in any way (even approximately):
there are many instances where \eqref{eq:obj_bm} has a benign landscape that is not explained by our results.
For example, much existing work (e.g., \cite{Bandeira2016a}, or the vast literature on Kuramoto oscillators) concerns the case $r = 2$, which our result does not cover.
This is particularly relevant to our results' robustness to monotone adversaries,
as there are problem instances (not covered by our results) where the nonconvex approach succeeds but is \emph{not} robust to monotone adversaries; see \Cref{sec:monotone}.
On the other hand, this illustrates the benefits of requiring $r \geq 3$: the Kuramoto oscillator literature presents many examples (even with $\Delta = 0$) for which \eqref{eq:obj_bm} does \emph{not} have a benign landscape for $r = 2$,
and the known examples which are not robust to monotone adversaries also have $r = 2$.

The ``benign landscape'' aspect of our results (i.e., that second-order critical points are global optima) depends fundamentally on the fact that we can recover the ground truth \emph{exactly}.
Without this, landscape analysis becomes more difficult, as it is more complicated to describe the global optimum.
This can be seen in the literature on more general orthogonal group synchronization, where we do not expect exact recovery in the presence of noise and, consequently, the best-known conditions for a benign landscape are more strict (see, e.g., the non-exact--recovery results of \citet{McRae2024,Ling2023b}).
Nevertheless, outside the exact recovery regime, one could obtain weaker statistical recovery results (second-order critical points are close to the ground truth)
as in, for example, \citet{Bandeira2016a,Ling2023a,McRae2024}, but we do not explore this here.

We now turn to the proof of \Cref{thm:z2_determ}.
The proof resembles that of \citet{McRae2024}, but, because we are optimizing over spheres rather than more general Stiefel manifolds,
we can refine that proof to obtain a stronger result.

The first step is to write down the conditions a second-order critical point $Y$ of \eqref{eq:obj_bm} must satisfy.
See \cite{Boumal2019} for an overview for problems of this form.
In short, \eqref{eq:obj_bm} is an optimization problem over a Riemannian manifold (in particular, a product of $n$ $(r-1)$-spheres embedded in $\R^r$).
Setting
\[
	S(Y) \coloneqq \ddiag(C Y Y^\top) - C,
\]
a second-order critical point $Y$ must satisfy the following:
\begin{itemize}
	\item $S(Y) Y = 0$; this is equivalent to the Riemannian gradient being zero.
	\item For any $\Ydt \in \R^{n \times r}$ satisfying $\diag(Y \Ydt^\top) = 0$, we have $\ip{S(Y)}{\Ydt \Ydt^\top} \geq 0$;
	this is equivalent to the Riemannian Hessian being negative semidefinite.
\end{itemize}

We can, without loss of generality (see \cite{Bandeira2016a,McRae2024} for similar arguments), assume that $z=\onevec$, as the change of variables $Y \mapsto \diag(z) Y$ implies the result for a general $z\in \{\pm 1\}^n$.
Although $\Delta$ changes, the important quantities $\opnorm{\Delta}$ and $\rho^\Delta$ do not.
Under this assumption, we have $C = A + \Delta$, and $\rho^\Delta_i = - \sum_{j = 1}^n \Delta_{ij}$ (recalling that we have assumed $\Delta$ is zero on the diagonal).

Let $L^\Delta = \diag(\Delta \onevec) - \Delta$ be the Laplacian matrix\footnote{This Laplacian matrix is not, in general, positive semidefinite like a graph Laplacian, as the entries of $\Delta$ are not assumed to be nonnegative.} of $\Delta$, and note that $L^\Delta \onevec = 0$ and
\begin{equation*}
	L^\Delta_{ij} = \begin{cases} -\rho^\Delta_i & \text{if } i = j, \\
		-\Delta_{ij} & \text{otherwise}.
	\end{cases}
\end{equation*}
Note that for $i \neq j$, $C_{ij} = - (L_{ij} + L^\Delta_{ij})$.
Furthermore, note that adding diagonal entries to $C$ has no effect on $S(Y)$ due to the constraint $\ddiag(Y Y^\top) = I_n$.
We can therefore write
\[
	S(Y) = \Lhat - \ddiag(\Lhat Y Y^\top),
\]
where $\Lhat := L + L^\Delta$ .

We will use the notation
\[
Y = \begin{bmatrix*}
	Y_1 \\ \vdots \\ Y_n
\end{bmatrix*}
\in \R^{n\times r}\qquad \text{and} \qquad
\Ydt = \begin{bmatrix*}
	\Ydt_1 \\ \vdots \\ \Ydt_n
\end{bmatrix*}
\in \R^{n\times r}.
\]
We will also treat the $1 \times r$ matrix rows as vectors in $\R^r$.

To prove Theorem \ref{thm:z2_determ}, it suffices to show that $Y=\onevec u^\top$ for some $u\in \R^r$;
the problem constraint then ensures $\norm{u}^2 = Y_1 Y_1^\top = 1$.
To show this, our proof will consider a general decomposition $Y=\onevec u^\top + W$, where $u \in \R^r$ is chosen such that $W^\top \onevec = 0$. The desired result is thus equivalent to $W = 0$.

With this decomposition defined, we can state the following technical lemma (proved at the end of this section) that we will need:
\begin{lemma}
	\label{lem:Qnorm_bd}
	Given a feasible point $Y$ of \eqref{eq:obj_bm},
	let $Q$ be the $n \times n$ matrix with entries
	\[
	Q_{ij}=\frac{1}{4}\norm{Y_i-Y_j}^4.
	\]
	Decomposing $Y = \onevec u^\top + W$ with $W^\top \onevec = 0$
	and setting $a_i = \frac{1}{4}\norm{W_i}^4$ (where $W_i \in \R^r$ is the $i$th row of $W$) for $i = 1, \dots, n$, we have
	\[
	Q = a \onevec^\top + \onevec a^\top + \Qtl,
	\]
	where $\Qtl$ is a matrix satisfying $\nucnorm{\Qtl} \leq 14 \normF{W}^2$.
\end{lemma}

The key step in the proof is to apply the second-order criticality inequality $\ip{S(Y)}{\Ydt \Ydt^\top} \geq 0$
for a particular choice of $\Ydt \in \R^{n \times r}$ satisfying the requirement $\diag(Y \Ydt^\top) = 0$.
It is, in fact, convenient to choose $\Ydt$ \emph{randomly} and take an expectation.
Our random choice of $\Ydt$ is also used by \citet{Ling2023b} and is a special case of those of \citet{McRae2024,Ling2023a}.
We choose
\[
\Ydt_i = \Gamma - \ip{\Gamma}{Y_i} Y_i,
\]
where $\Gamma$ is a $1 \times r$ matrix whose entries are i.i.d.\ standard Gaussian random variables. Note that $\ip{Y_i}{\Ydt_i} = 0$ for all $i$, so indeed we have $\diag(Y \Ydt^\top) = 0$.

Next, it follows directly from the choice of $\Ydt$ that
\begin{align*}
	\E \ip{\Ydt_i}{\Ydt_j}
	&= r - 2 + \ip{Y_i}{Y_j}^2 \\
	&= r - 2+ \parens*{1 - \frac{1}{2} \norm{Y_i - Y_j}^2}^2 \\
	&= r - 1 - \norm{Y_i - Y_j}^2 + \underbrace{\frac{1}{4} \norm{Y_i - Y_j}^4}_{Q_{ij}} \\
	&= r - 3 + 2 \ip{Y_i}{Y_j} + Q_{ij}.
\end{align*}
Equivalently, $\E \Ydt \Ydt^\top = (r - 3) \onevec \onevec^\top + 2 Y Y^\top + Q$. We apply the second-order criticality condition under the expectation to obtain 
\begin{align*}
	0 &\leq \ip{S(Y)}{\E \Ydt \Ydt^\top} \\
	&= \ip{S(Y)}{(r-3) \onevec \onevec^\top + 2 Y Y^\top + Q} \\
	&\overset{\mathclap{\text{(a)}}}{=} \ip{S(Y)}{ (r-3) \onevec \onevec^\top  + Q} \\
	&= \ip{ L + L^\Delta - \ddiag((L + L^\Delta) Y Y^\top)}{ (r-3) \onevec \onevec^\top  + Q } \\
	&\overset{\mathclap{\text{(b)}}}{=} -\ip{ A + \Delta}{Q} - (r-3) \ip{L + L^\Delta}{Y Y^\top}.
\end{align*}
Equality (a) uses the fact that $\ip{S(Y)}{Y Y^\top} = 0$ (this is clear from the first-order criticality condition $S(Y) Y = 0$, but one can easily verify that the same holds for any feasible $Y$).
Equality (b) uses the facts that $\onevec$ is in the null space of both $L$ and $L^\Delta$ and that $Q$ is zero on the diagonal. Rearranging, we obtain the inequality
\begin{equation}
	\label{eq:randYdt_basic_ineq}
	(r - 3) \ip{L + L^\Delta}{Y Y^\top} + \ip{A}{Q} \leq -\ip{ \Delta }{Q}.
\end{equation}
Recalling the decomposition $Y = \onevec u^\top + W$ with $W^\top \onevec = 0$,
note that
\begin{align}
	\label{eq:LYY_lb}
	\ip{L + L^\Delta}{Y Y^\top}
	&= \ip{L + L^\Delta}{W W^\top} \nonumber \\
	&= \ip{L - \Delta}{W W^\top} - \sum_{i=1}^n \rho^\Delta_i \norm{W_i}^2 \nonumber \\
	&\geq \parens{ \lambda_2 - \opnorm{\Delta}} \normF{W}^2 - \sum_{i=1}^n \rho^\Delta_i \norm{W_i}^2. 
\end{align}
The last inequality follow from the fact that for any symmetric matrices $A, B$ with $B \succeq 0$, $\ip{A}{B} \geq \lambda_{\min}(A) \tr(B)$.
Because $G$ is a connected graph, unless $Y_1 = \cdots = Y_n$ (or, equivalently, $W = 0$), we have
\begin{equation}
	\label{eq:quartic_neg}
	\ip{A}{Q} = \frac{1}{4} \sum_{i,j=1}^n A_{ij} \norm{Y_i - Y_j}^4 > 0.
\end{equation}
From now on, assume, by way of contradiction, that $W \neq 0$, so \eqref{eq:quartic_neg} holds.

Next, we upper bound $-\ip{\Delta}{Q}$. By \Cref{lem:Qnorm_bd},
\begin{align}
	\label{eq:deltaQ_bd}
	\ip{-\Delta}{Q}
	&= -2 \ip{\Delta}{a \onevec^\top} - \ip{\Delta}{\Qtl} \nonumber \\
	&\leq -2 \sum_{i,j=1}^{n} \Delta_{ij} a_i + 14 \opnorm{\Delta} \normF{W}^2 \nonumber \\
	&=  \frac{1}{2} \sum_{i=1}^n \rho^\Delta_i \norm{W_i}^4 + 14 \opnorm{\Delta} \normF{W}^2. 
\end{align}
Combining \eqref{eq:randYdt_basic_ineq}, \eqref{eq:LYY_lb}, \eqref{eq:quartic_neg}, and \eqref{eq:deltaQ_bd}, we obtain
\begin{align*}
	(r-3) (\lambda_2 - \opnorm{\Delta}) \normF{W}^2
	&<14 \opnorm{\Delta} \normF{W}^2 + \sum_{i=1}^n \rho^\Delta_i \brackets*{ (r-3) \norm{W_i}^2 + \frac{1}{2} \norm{W_i}^4 } \\
	&\leq 14 \opnorm{\Delta} \normF{W}^2 + \rho^\Delta \brackets*{ (r-3) \normF{W}^2 + \frac{1}{2} \sum_{i=1}^n\norm{W_i}^4 }.
\end{align*}
Using the fact that $\norm{W_i} \leq 2$ and therefore $\norm{W_i}^4 \leq 4 \norm{W_i}^2$, we obtain
\begin{align*}
	(r - 3)(\lambda_2 - \rho^\Delta - \opnorm{\Delta}) \normF{W}^2
	< 14 \opnorm{\Delta} \normF{W}^2 + 2 (\rho^\Delta \vee 0) \normF{W}^2,
\end{align*}
where we have used the shorthand $a \vee b = \max\{a, b\}$.
The last inequality is equivalent to
\begin{equation*}
	\brackets*{ (r-3) [\lambda_2 - \rho^\Delta] - 2(\rho^\Delta \vee 0) - (r + 11) \opnorm{\Delta} }\normF{W}^2 < 0.
\end{equation*}
If
\begin{equation}
	\label{eq:z2_cond_full}
	(r-3) [\lambda_2 - \rho^\Delta] - 2(\rho^\Delta \vee 0) - (r + 11) \opnorm{\Delta} \geq 0,
\end{equation}
we obtain a contradiction (of the supposition $W \neq 0$) and thus must have $W = 0$.
It is clear that \eqref{eq:z2_cond_full} holds if $r = 3$ and $\Delta = 0$ (and hence $\rho^\Delta = 0$).
If $r \geq 4$, considering the two cases $\rho^\Delta \geq 0$ and $\rho^\Delta \leq 0$ reveals that \eqref{eq:z2_cond} implies \eqref{eq:z2_cond_full}.
This completes the proof of \Cref{thm:z2_determ}.

Finally, we prove the technical lemma we used:
\begin{proof}[Proof of \Cref{lem:Qnorm_bd}]
	Expanding the entries of $Q$, we have
	\begin{align*}
		Q_{ij}
		&= \frac{1}{4} \norm{Y_i - Y_j}^4 \\
		&= \frac{1}{4} \norm{W_i - W_j}^4 \\
		&= \frac{1}{4} \parens*{ \norm{W_i}^2 + \norm{W_j}^2 - 2 \ip{W_i}{W_j} }^2 \\
		&= \frac{1}{4} \parens*{ \norm{W_i}^4 + \norm{W_j}^4 + 2 \norm{W_i}^2 \norm{W_j}^2 + 4 \ip{W_i}{W_j}^2 - 4(\norm{W_i}^2 + \norm{W_j}^2) \ip{W_i}{W_j} }\\
		&= (a \onevec^\top + \onevec a^\top)_{ij} + \underbrace{\parens*{ \frac{1}{2} \norm{W_i}^2 \norm{W_j}^2 +  \ip{W_i}{W_j}^2 - (\norm{W_i}^2 + \norm{W_j}^2) \ip{W_i}{W_j}}}_{\eqqcolon \Qtl_{ij}}.
	\end{align*}
	Now we must bound $\nucnorm{\Qtl}$.
	First, noting that for a matrix $M \succeq 0$, $\nucnorm{M} = \tr(M)$,
	we have (again recalling that each $\norm{W_i} \leq 2$)
	\begin{align*}
		\nucnorm*{ \brackets*{ \frac{1}{2} \norm{W_i}^2 \norm{W_j}^2 + \ip{W_i}{W_j}^2 }_{ij} }
		&= \frac{3}{2}\sum_{i=1}^n \norm{W_i}^4 \\
		&\leq 6\normF{W}^2.
	\end{align*}
	To bound the remaining terms $\nucnorm{ \brackets{ \norm{W_i}^2 \ip{W_i}{W_j} }_{ij} } = \nucnorm{ \brackets{ \norm{W_j}^2 \ip{W_i}{W_j} }_{ij} }$,
	note that the matrix $\brackets{\norm{W_i}^2 \ip{W_i}{W_j} }_{ij}$ is the Hadamard (entrywise) product of $W W^\top$ (whose nuclear norm is $\normF{W}^2$) with the matrix $v \onevec^\top$, where $v_i = \norm{W_i}^2$.
	The inequality \citep[Theorem 5.6.2]{Horn1991} bounds the singular values of Hadamard products with such structure: it implies that, because $\abs{v_i} \leq 4$ for all $i$,
	\begin{align*}
		\nucnorm*{ \brackets*{ \norm{W_i}^2 \ip{W_i}{W_j} }_{ij} }\leq 4 \normF{W}^2.
	\end{align*}
	Putting these bounds together with the triangle inequality, we obtain the bound $\nucnorm{\Qtl} \leq 14 \normF{W}^2$, completing the proof of \Cref{lem:Qnorm_bd}.
\end{proof}

\section{Asymptotic probabilistic results for $\Z_2$ synchronization}
\label{sec:z2_asymp}
In this section, we apply the deterministic result \Cref{thm:z2_determ} to the problems with random noise and graphs described in \Cref{sec:intro_gauss,sec:intro_er}.
\subsection{Complete graph and Gaussian noise}
\label{sec:z2_proof_gauss}
First, we apply \Cref{thm:z2_determ} to the model \eqref{eq:z2_gauss_model} in \Cref{sec:intro_gauss} to prove \Cref{cor:z2_gauss};
this is not only a simple and illustrative application of \Cref{thm:z2_determ} but also a warm-up to the more involved probabilistic analyses later in this paper.
In this case, the measurement graph $G$ is the complete graph on the vertices $1, \dots, n$,
and the noise matrix has the form $\Delta = \sigma W$, where $W$ is an $n \times n$ symmetric random matrix with i.i.d.\ (modulo symmetry) standard normal entries.
To apply \Cref{thm:z2_determ}, we need to estimate three quantities: $\lambda_2$, $\rho^{\Delta}$ and $\opnorm{\Delta}$.

Because $G$ is the complete graph, it has graph Laplacian $L = n I_n - \onevec \onevec^\top$ and thus algebraic connectivity $\lambda_2 = n$. Next, noting that $\rho^\Delta$ in \eqref{eq:rhomax} is the maximum of $n$ zero-mean Gaussian random variables with variance $(n-1) \sigma^2$,
standard concentration inequalities for Gaussian variables and matrices
(see, e.g., \cite[Section 3.2]{Bandeira2018})
yield that, for any $\epsilon' > 0$, with probability $\to 1$ as $n \to \infty$,
\begin{equation*}
	\rho^{\Delta} \leq \sigma\sqrt{(2+\epsilon')n\log n}\quad \text{and} \quad \opnorm{\Delta} \leq  3\sigma \sqrt{n}.
\end{equation*}
On this event, plugging $\lambda_2=n$ and these upper bounds for $\rho^\Delta$ and $\opnorm{\Delta}$ into \Cref{thm:z2_determ}, we obtain that for $r\geq 4$, the desired result holds if
\begin{equation*}
	\sigma\sqrt{(2 + \epsilon')n\log n} + 45 \sigma \sqrt{n}  \leq \frac{r-3}{r-1}n.
\end{equation*}
The first term dominates as $n \to \infty$.
With, for example, $\epsilon' = \frac{\epsilon}{2}$ (with $\epsilon$ as in the corollary statement),
this holds for sufficiently large $n$ when $\sigma \leq \frac{r-3}{r-1} \sqrt{\frac{n}{(2+\epsilon)\log n}}$.
In the case $r = 3$, the condition implies $\sigma = 0$ and therefore $\Delta = 0$, so \Cref{thm:z2_determ} still applies.
This completes the proof of \Cref{cor:z2_gauss}.

\subsection{\ER{} random graph and Bernoulli noise}
\label{sec:z2_bern}
We now state and prove a general result for \ER{} random graphs with Bernoulli noise, as described in \Cref{sec:intro_er}:
\begin{theorem}
	\label{thm:z2_er_bern}
	Consider the $\Z_2$ synchronization measurement model with Bernoulli noise \eqref{eq:z2_bernoulli_model} for some $\delta \in (0, 1]$
	on an \ER{} random graph $G \sim \ERmath(n, p)$.

	Suppose that the problem parameters $p$, $\delta$, and $r \geq 4$ (which can depend on $n$) satisfy, for some $0 < \epsilon \leq 1/3$ independent of $n$,
	\begin{equation}
		\label{eq:z2_bern_cond}
		\frac{np}{\log n} \parens*{ 1 - \sqrt{1 - \parens*{ \frac{r-3}{r-1} - \epsilon}^2 \delta^2 } } \geq 1.
	\end{equation}
	Then, with probability $\to 1$ as $n \to \infty$,
	for cost matrix $C$ as in \eqref{eq:bern_er_C},
	every second-order critical point $Y$ of \eqref{eq:obj_bm} satisfies $Y = z u^\top$ for some unit-norm $u \in \R^r$.

	A more restrictive condition for which the same result holds is that, for some $\epsilon > 0$ independent of $n$,
	\begin{equation}
		\label{eq:z2_bern_cond_simple}
		\frac{1}{\delta} \leq \frac{r-3}{r-1} \sqrt{\frac{np}{(2+\epsilon) \log n}}.
	\end{equation}
\end{theorem}
Note that condition \eqref{eq:z2_bern_cond} is not tight in terms of $r$ in the ``noiseless'' case $\delta = 1$,
as then \Cref{thm:z2_determ} guarantees that even $r=3$ suffices when the measurement graph is connected (which occurs when $\frac{np}{\log n} > 1$).
In fact, by quite different methods, \citet{Abdalla2022} show that even $r = 2$ suffices under similar conditions on $p$ and $n$.

The threshold $\frac{np}{\log n} (1 - \sqrt{1 - \delta^2}) = 1$,
which \Cref{thm:z2_er_bern} approaches for large $r$ and small $\epsilon$,
is the best possible for exact recovery \cite{Hajek2016a,Jog2015}.
In particular, \citet{Hajek2016a} showed that (in the case that $p$ is proportional to $\frac{\log n}{n}$) this threshold is achieved by the semidefinite relaxation \eqref{eq:obj_sdp};
we show that we can reach the same threshold via nonconvex optimization of \eqref{eq:obj_bm}.
Furthermore, our result, like the SDP approach, is robust to monotone adversaries: see \Cref{sec:monotone} for more details.

To our knowledge, the only previous work to consider our nonconvex optimization approach for this specific model is that of \citet[Cor.~2.3]{Ling2023b}; though phrased differently in terms of Kuramoto oscillators, this is equivalent to our model in the specific case $p = 1$ (complete graph).
Our result improves this previous result by a constant factor.
%

We now prove \Cref{thm:z2_er_bern} from \Cref{thm:z2_determ}.
As in the proof of \Cref{thm:z2_determ} in \Cref{sec:z2_deterministic} we can assume, without loss of generality, that $z = \onevec$.
If $A$ is the adjacency matrix of the \ER{} random graph $G$,
the model \eqref{eq:z2_bernoulli_model} with cost matrix $C$ as in \eqref{eq:bern_er_C} gives, independently (modulo symmetry) for all $i$ and $j$,
\[
	C_{ij} = \begin{cases}
		A_{ij} & \text{with probability } \frac{1 + \delta}{2} \\
		- A_{ij} & \text{with probability } \frac{1 - \delta}{2}.
	\end{cases}
\]
We could directly apply \Cref{thm:z2_determ} with adjacency matrix $A$ and noise matrix $\Delta = C - A$.
However, this turns out to be problematic for two reasons:
first, this $\Delta$ is not zero-mean, and thus $\opnorm{\Delta}$ may be quite large (though this problem could be remedied with some rescaling tricks similar to what we do below);
second and more seriously, the tight coupling between the graph $G$ and the locations of the errors turns out to make the general result \Cref{thm:z2_determ} not as tight as possible when applied directly.

We therefore take a different approach.
To apply \Cref{thm:z2_determ}, we can decompose $C$ as $C = \Atl + \Deltatl$ however we like as long as the entries of $\Atl$ are nonnegative.
The adjacency matrix $A$ of an \ER{} graph is well-approximated (up to rescaling) by that of the complete graph.
In particular, noting that $\E C = \delta p (\onevec \onevec^\top - I_n)$,
we write
\begin{align}
	C &= \underbrace{\delta p (\onevec \onevec^\top - I_n)}_{\Atl} + \underbrace{C - \E C}_{\Deltatl}. \label{eq:Cdecomp_bern}
\end{align}
Now, $\Atl$ is the (rescaled) complete-graph adjacency matrix,
and $\Deltatl$ includes both the original measurement error and the ``sampling noise'' of the random graph itself.
The scaled graph Laplacian of $\Atl$ is $\Ltl = \delta p (n I_n - \onevec \onevec^\top)$,
for which $\lambda_2 = \cdots = \lambda_n = n \delta p$.
We apply \Cref{thm:z2_determ} with this $\Atl$ and $\Deltatl$.

We first need to bound the error operator norm $\opnorm{\Deltatl}$.
The following general-purpose result is useful:
\begin{lemma}
	\label{lem:X_conc}
	Let $X$ be a real symmetric $n \times n$ random matrix with independent (modulo symmetry) and zero-mean entries $X_{ij}$ satisfying $\abs{X_{ij}} \leq 2$ almost surely and $\E X_{ij}^2 \leq v$ for all $i,j$,
	where $v \geq \frac{\log n}{n}$.
	Then, with probability at least $1 - n^{-3}$,
	\[
	\opnorm{X} \lesssim \sqrt{nv}.
	\]
\end{lemma}
This is a slight generalization of \cite[Theorem 5]{Hajek2016}
and is a corollary of, for example, \cite[Lemma 4.1]{Bandeira2018}.
Using the inequality $1 - \sqrt{1 - x} \leq x$ for $x \in [0, 1]$,
condition \eqref{eq:z2_bern_cond} implies $\frac{np}{\log n} \delta^2 \geq 1$;
this also implies $p \geq \frac{\log n}{n}$.
We can then apply \Cref{lem:X_conc} with $X = \Deltatl = C - \E C$ and $v = p$;
it is clear from the definition of $\Deltatl$ that $\E \Deltatl_{ij} = 0$ and $\abs{\Deltatl_{ij}} \leq 2$, while
\[
	\E \Deltatl_{ij}^2 \leq \E C_{ij}^2 = \E A_{ij} = p.
\]
We then obtain, with probability $\to 1$ in $n$,
\begin{align*}
	\opnorm{\Deltatl} \lesssim \sqrt{np} \leq \frac{n \delta p}{\sqrt{\log n}}.
\end{align*}
This implies that for large enough $n$, the $\opnorm{\Deltatl}$ term in the condition \eqref{eq:z2_cond} is negligible compared to $\lambda_2 = n\delta p$.
It thus remains to deal with $\rho^{\Deltatl}$ from \eqref{eq:z2_cond}.
It suffices to show that, with probability $\to 1$,
\[
	\frac{\rho^{\Deltatl}}{ n\delta p} \leq \frac{r-3}{r-1} - \epsilon'
\]
for any $\epsilon' > 0$ independent of $n$.
This $\epsilon'$ provides some slack to account for the asymptotically negligible $\opnorm{\Deltatl}$.
It turns out we also need slack in another place in the proof, so we set, with $\epsilon > 0$ as in the theorem statement,
\begin{equation}
	\label{eq:z2_proof_cs}
	c \coloneqq \frac{r-3}{r-1} - \epsilon, \qquad c' \coloneqq \frac{r-3}{r-1} - \frac{\epsilon}{2},
\end{equation}
and we will show that $\rho^{\Deltatl} \leq c' n\delta p$ with high probability.

To do so, we need to upper bound the $n$ random variables
\[
	\rho^{\Deltatl}_i = -\sum_{j = 1}^n (C_{ij} - \E C_{ij}) = \sum_{\substack{j = 1 \\ j \neq i}}^n (\delta p - C_{ij}).
\]
For each $i$, the sum in $j$ is over $n-1$ i.i.d.\ random variables.
We can calculate, for $t \in \R$,
\begin{align*}
	\E e^{t(\delta p - C_{ij})}
	&= e^{t \delta p} \brackets*{ p \frac{1 + \delta}{2} e^{-t} + p \frac{1 - \delta}{2} e^t + 1 - p},
\end{align*}
so, using the inequality $\log (1 + x) \leq x$,
\begin{align*}
	\log(\E e^{t \rho^{\Deltatl}_i})
	&= (n-1) \log \parens*{ e^{t \delta p} \brackets*{ p \frac{1 + \delta}{2} e^{-t} + p \frac{1 - \delta}{2} e^t + 1 - p} } \\
	&\leq np \brackets*{ \delta t + \frac{1 + \delta}{2} e^{-t} + \frac{1 - \delta}{2} e^t - 1  }.
\end{align*}
Recalling from \eqref{eq:z2_proof_cs} that $c = c' - \frac{\epsilon}{2}$, we then obtain, by a Chernoff bound and union bound, for $t \geq 0$,
\begin{align*}
	\P(\rho^{\Deltatl} \geq c' n\delta p)
	&\leq n \frac{ \max_i~\E e^{t \rho^{\Deltatl}_i} }{ e^{c' n\delta p t} } \\
	&\leq n \exp \parens*{ np \brackets*{ \delta (1 - c') t + \frac{1 + \delta}{2} e^{-t} + \frac{1 - \delta}{2} e^t - 1  } } \\
	&= \exp\parens*{ \log n + np \brackets*{ (1 - c) \delta t + \frac{1 + \delta}{2} e^{-t} + \frac{1 - \delta}{2} e^t - \frac{\epsilon \delta t}{2} - 1 } } \\
	&\leq \exp\parens*{ \log n + np \brackets*{ \frac{1 + c\delta}{2} e^{-t} + \frac{1 - c\delta}{2} e^t - \frac{\epsilon \delta t}{2} - 1} },
\end{align*}
where the last inequality uses the fact that, for all $t \geq 0$, $t \leq \sinh t = \frac{e^t - e^{-t}}{2}$.

Choosing $t = \frac{1}{2} \log \frac{1 + c\delta}{1 - c \delta}$,
we obtain
\begin{align*}
	\P(\rho^{\Deltatl} \geq c' n\delta p)
	&\leq \exp\parens*{ \log n + np \brackets*{ \sqrt{1 - c^2 \delta^2} - 1} - \frac{\epsilon \delta np}{4} \log \frac{1 + c \delta}{1 - c\delta}} \\
	&\leq \exp\parens*{ - \frac{\epsilon c \delta^2 np}{2} }.
\end{align*}
The second inequality uses $\frac{1}{2} \log \frac{1 + x}{1-x} \geq x$ for $x \geq 0$ and the condition \eqref{eq:z2_bern_cond}.
Noting that the same condition \eqref{eq:z2_bern_cond} implies
\[
	\log n \leq np \parens*{ 1 - \sqrt{1 - c^2 \delta^2}} \leq np c^2 \delta^2
\]
(using the inequality $1 - \sqrt{1 - x} \leq x$ for $x \in [0,1]$)
and that $c \leq 1$, we obtain
\[
	\P(\rho^{\Deltatl} \geq c' n\delta p)
	\leq e^{-\frac{\epsilon}{2} \log n} \to 0.
\]
This proves the sufficiency of condition \eqref{eq:z2_bern_cond} in \Cref{thm:z2_er_bern}.
For the simplified condition \eqref{eq:z2_bern_cond_simple}, note that
\[
	\frac{np}{2 \log n} c^2 \delta^2 \geq 1
\]
implies \eqref{eq:z2_bern_cond} by the inequality $1 - \sqrt{1 - x} \geq \frac{x}{2}$ for $x \in [0, 1]$.
Redefining $\epsilon > 0$ as needed, we see that condition \eqref{eq:z2_bern_cond_simple} is sufficient.

\section{Graph clustering under the SBM}
\label{sec:sbm_main}
In this section, we adapt our arguments and results in \Cref{sec:z2_deterministic,sec:z2_asymp} to the problem of graph clustering under the binary stochastic block model (SBM) described in \Cref{sec:intro_sbm}.
The main result, of which \Cref{cor:z2_bern} is a simplified version, is the following theorem:
\begin{theorem}
	\label{thm:sbm_asymp}
	Consider the model \eqref{eq:sbm_model} with parameters $0 \leq q < p \leq 1$.
	Suppose the problem parameters $p$, $q$, and $r \geq 4$ (which can depend on the problem size $n$) satisfy, for some $0 < \epsilon \leq 1/12$ independent of $n$,
	\begin{equation}
		\label{eq:sbm_asymp_cond}
		\frac{n}{\log n} \parens*{ \sqrt{ p - \parens*{\frac{1}{r-1} + \epsilon} (p - q) } - \sqrt{q + \parens*{\frac{1}{r-1} + \epsilon} (p - q) } }^2 \geq 2.
	\end{equation}
	Then, with probability $\to 1$ as $n \to \infty$, with cost matrix $C$ as in \eqref{eq:sbm_C} or \eqref{eq:sbm_C_pqknown}, every second-order critical point $Y$ of \eqref{eq:obj_bm} satisfies $Y = z u^\top$ for some unit-norm $u \in \R^r$.
\end{theorem}
For sufficiently large $r$ and small $\epsilon$,
this approaches the optimal threshold $\frac{n}{\log n} (\sqrt{p} - \sqrt{q})^2 = 2$ for exact recovery \cite{Mossel2015,Abbe2016}.
The SDP relaxation \eqref{eq:obj_sdp}, which has the key benefit of robustness to monotone adversaries (see \Cref{sec:monotone}),
was shown to achieve this threshold by \citet{Hajek2016,Bandeira2018}.
Our result \Cref{thm:sbm_asymp} shows that we can achieve the same threshold (with similar robustness to monotone adversaries---see \Cref{thm:monotone} in \Cref{sec:monotone}) with continuous, benignly nonconvex optimization
while optimizing over fewer variables.

The first theoretical work on the nonconvex approach to this problem was by \citet{Bandeira2016a}, who studied the case $r = 2$.
That paper gives sufficient conditions for high-probability exact recovery, but these conditions scale suboptimally in the problem size $n$.
Recently, \citet{Ling2023b} proved the optimal dimension scaling,
showing that, for $r \geq 4$ and for some $c > 0$, $\frac{n}{\log n} \frac{(p-q)^2}{p+q} \geq c$ suffices for exact recovery.
Our result closes the remaining gap to the optimal threshold.

There has been a variety of work on more general graph clustering models that encompass both the $\Z_2$ synchronization problem with Bernoulli noise of \Cref{sec:intro_er} and the ordinary clustering problem of \Cref{sec:intro_sbm} as special cases (e.g., the ``weighted'' SBM of \citet{Jog2015}, the ``labeled'' SBM of \citet{Yun2016}, or the ``signed'' SBM of \citet{Wang2022}).
However, to combine or generalize these models in our algorithmic framework, we must reconcile two factors:
\begin{itemize}
	\item In the Bernoulli noise model, the problem difficulty is invariant to the actual value of the ground truth labels $z$ (this is clear from the proofs, in which we assume without loss of generality that $z = \onevec$).
	\item For the ordinary SBM, the cluster sizes matter greatly.
	It is standard in the literature to assume, as we have, that the clusters have equal size (indeed, this assumption is present even in \cite{Jog2015,Wang2022}).
	Otherwise, the problem difficulty depends on the relative cluster sizes (see, e.g., \cite{Hajek2016a,Yun2016}). Furthermore, if the parameters $p$ and $q$ are unknown, it becomes necessary to estimate the cluster sizes, which adds considerable technical complication to the algorithm and analysis (again, see \cite{Hajek2016a}).
\end{itemize}
Thus to unify the problems, we need either to impose a balancing requirement (which is artificial and unnecessary for the pure synchronization/``signed clustering'' problem) or to deal with the complications of the SBM with unbalanced clusters.
For simplicity and clarity, we avoid this in the present work.

\subsection{Deterministic condition for cluster recovery}
As the first step of proving \Cref{thm:sbm_asymp}, we consider how to translate the deterministic analysis of \Cref{sec:z2_deterministic} to the somewhat different graph clustering problem.
Recall from \Cref{sec:intro_sbm} that we are given the $n \times n$ adjacency matrix $A$ of a random graph $G$ with distribution given by \eqref{eq:sbm_model} for unknown cluster labels $z \in \{ \pm 1 \}^n$
and with (potentially unknown) problem parameters $0 \leq q < p \leq 1$.
We have further assumed that the clusters are balanced: $\ip{z}{\onevec} = 0$ (hence $n$ must be even).
We then attempt to recover $z$ by the problem \eqref{eq:obj_bm} with cost matrix $C = A - \frac{1}{n^2} \ip{A}{\onevec \onevec^\top} \onevec \onevec^\top$ as in \eqref{eq:sbm_C} or, if the parameters are known, $C = A - \frac{p+q}{2} \onevec \onevec^\top$ as in \eqref{eq:sbm_C_pqknown}.

We can find a condition similar to the local stability for $\Z_2$ synchronization in \eqref{eq:z2_stability_cond} as follows:
for every $i\in \{1, \dots, n\}$, consider the quantity 
\begin{align*}
	d^z_i
	\coloneqq&~\abs{\{ \text{neighbors of $i$ from \textbf{same} cluster} \}} - \abs{\{ \text{neighbors of $i$ from \textbf{other} cluster} \}} \\
	=&~z_i \sum_{\substack{j = 1\\j \neq i}}^n A_{ij} z_j.
\end{align*}
A necessary condition for $\pm z$ to be the unique optimum of \eqref{eq:obj_orig} is that flipping the sign of one entry of $z$ must strictly decrease the objective function.
One can verify that, for any $C$ of the form $C = A - \alpha \onevec \onevec^\top$,
this is equivalent to $d^z_i + 2 \alpha > 0$ for all $i$.
Thus we expect $d^z_{\min}:=\min_i~d^z_i$ to play a role in the analysis.
These quantities have appeared before in the literature. One result of \citet{Mossel2015} says that we can asymptotically obtain exact recovery if and only if $n \P(d^z_i \leq 0) \to 0$ as $n \to \infty$.
For the SDP relaxation, \citet{Hajek2016} show that $d^z_{\min} > \opnorm{A - \E A}$ is a (deterministically) sufficient condition for exact recovery.
Our first step is to show that a similar condition suffices for the nonconvex problem \eqref{eq:obj_bm}.
A careful application of \Cref{thm:z2_determ} yields the following result:
\begin{theorem}
	\label{cor:sbm_determ}
	Let $r\geq 4$ be an integer, suppose $p > q$, and suppose one of the following holds:
	\begin{itemize}
		\item We set the cost matrix $C$ as in \eqref{eq:sbm_C} and assume
		\begin{equation}
			\label{eq:sbm_cond}
			d^z_{\min} \geq \frac{r + 11}{r-3} (2 \opnorm{A - \E A} + p) + \frac{n(p-q)}{r-1}.
		\end{equation}
		\item We set $C$ as in \eqref{eq:sbm_C_pqknown} and assume
		\begin{equation}
			\label{eq:sbm_cond_pqknown}
			d^z_{\min} \geq \frac{r + 11}{r-3} \opnorm{A - \E A} + \frac{n(p-q)}{r-1}.
		\end{equation}
	\end{itemize}
	Then every second-order critical point $Y$ of \eqref{eq:obj_bm} satisfies $Y = z u^\top$ for some unit-norm $u \in \R^r$.
	
\end{theorem}
In the limit $r \to \infty$ (which one can take even for fixed $n$ for the full SDP relaxation---see, e.g., \cite{McRae2024}), \eqref{eq:sbm_cond_pqknown} becomes $d^z_{\min} > \opnorm{A - \E A}$;
this is identical to the requirement of \citet{Hajek2016}.
The extraneous factor of two in \eqref{eq:sbm_cond} (the $p$ term is negligible for large $n$) is the result of a crude approximation and could be improved at the cost of a slightly more complicated theorem statement.
But in any case, the $\opnorm{A - \E A}$ terms of \eqref{eq:sbm_cond} and \eqref{eq:sbm_cond_pqknown} will prove to be negligible in our asymptotic analysis.

\begin{proof}[Proof of \Cref{cor:sbm_determ}]
The result follows from applying \Cref{thm:z2_determ} to a suitably translated problem.
Similarly to the proof of \Cref{thm:z2_er_bern} in \Cref{sec:z2_bern}, we compare the cost matrix $C$ to that arising from an idealized complete-graph problem.
We cannot assume now that $z = \onevec$ (as the balancing property $\ip{z}{\onevec} = 0$ is important), so the notation will be somewhat more cumbersome than before.

First, we prove the first part of \Cref{cor:sbm_determ}, where $C = A - \frac{1}{n^2} \ip{A}{\onevec \onevec^\top} \onevec \onevec^\top$.
For convenience, we denote
\[
\Pone \coloneqq \frac{1}{n} \onevec \onevec^\top \quad \text{and} \quad \Poneperp \coloneqq I_n - \Pone = I_n - \frac{1}{n} \onevec \onevec^\top
\]
as the orthogonal projection matrices onto $\spn\{\onevec\}$ and $\spn\{\onevec\}^\perp$ respectively.
We can then more compactly write $C = A  - \Pone A \Pone$.
Recalling the expression for $\E A$ in \eqref{eq:EA}, note that
\begin{equation}
	\label{eq:sbm_EC}
	\E C
	= \frac{p-q}{2} z z^\top - p I_n + p \Pone.
\end{equation}
We can then decompose
\begin{align*}
	C
	&= \frac{p-q}{2} z z^\top - p I_n + p \Pone + C - \E C \\
	&= \underbrace{\frac{p-q}{2} (z z^\top - I_n)}_{\diag(z) \Atl \diag(z) } - \frac{p+q}{2} I_n + \underbrace{p \Pone + C - \E C}_{\Deltatl},
\end{align*}
where $\Atl = \frac{p-q}{2} (\onevec \onevec^\top - I_n)$ is the complete-graph adjacency matrix with scaling factor $\frac{p-q}{2}$.
As before, we can assume the diagonal elements of $C$ are zero;
we therefore ignore the identity term in the decomposition of $C$ and any diagonal elements of $\Deltatl$.
We then apply \Cref{thm:z2_determ} with graph adjacency matrix $\Atl$ and noise matrix $\Deltatl$.

Clearly, $\lambda_2 = n\frac{p-q}{2}$.
Next, note that we can bound
\begin{align*}
	\opnorm{\Deltatl}
	= \opnorm{p \Pone + A - \E A - \Pone(A - \E A) \Pone}
	\leq p + 2 \opnorm{A - \E A}.
\end{align*}

Last, we need to bound $\rho^{\Deltatl} = \max_i~\rho^{\Deltatl}_i$ where
\[
\rho^{\Deltatl}_i = - z_i \sum_{\substack{j = 1\\j \neq i}}^n \Deltatl_{ij} z_j.
\]
Recalling \eqref{eq:sbm_EC}, note that for $i \neq j$,
\begin{align*}
	\Deltatl_{ij} &= \frac{p}{n} + C_{ij} - \E C_{ij} \\
	&= \frac{p}{n} + A_{ij} - \frac{1}{n^2} \sum_{k,\ell} A_{k\ell} - \frac{p-q}{2} z_i z_j - \frac{p}{n} \\
	&= A_{ij} - \frac{\ip{A}{\onevec \onevec^\top}}{n^2} - \frac{p-q}{2} z_i z_j.
\end{align*}
Then, noting that $\sum_{j\neq i} z_j = \ip{z}{\onevec} - z_i = -z_i$,
we can calculate
\begin{align*}
	\rho^{\Deltatl}_i
	&= -z_i \sum_{j \neq i} A_{ij} z_j + z_i\parens*{ \sum_{j\neq i} z_j }\frac{\ip{A}{\onevec \onevec^\top}}{n^2} + (n-1) \frac{p-q}{2} \\
	&= - d^z_i - \frac{\ip{A}{\onevec \onevec^\top}}{n^2} + (n-1) \frac{p-q}{2} \\
	&\leq - d^z_i + n \frac{p-q}{2}.
\end{align*}
Then, clearly, $\rho^{\Deltatl} \leq - d^z_{\min} + n \frac{p-q}{2}$.
To show the condition \eqref{eq:z2_cond} of \Cref{thm:z2_determ} holds,
it then suffices to have
\[
- d^z_{\min} + n \frac{p-q}{2} + \frac{r+11}{r-3} (p + 2 \opnorm{A - \E A}) \leq \frac{r-3}{r-1} n \frac{p-q}{2}.
\]
The condition \eqref{eq:sbm_cond} is equivalent.
This proves the sufficiency of \eqref{eq:sbm_cond} when $C$ is set by \eqref{eq:sbm_C_pqknown}.

The proof of the second part, where $C = A - \frac{p+q}{2} \onevec \onevec^\top$, is simpler.
Again recalling \eqref{eq:EA}, we now have
\[
\E C= \frac{p-q}{2} z z^\top - p I_n,
\]
so we can decompose
\begin{align*}
C &= \E C + C - \E C \\
&=\underbrace{\frac{p-q}{2} (z z^\top - I_n)}_{\diag(z) \Atl \diag(z) } - \frac{p+q}{2} I_n + \underbrace{A - \E A}_{\Deltatl}.
\end{align*}
Noting that now
\[
\rho^{\Deltatl}_i = - (d^z_i - \E d^z_i) = - d^z_i + n \frac{p-q}{2} - p,
\]
the result follows if
\[
- d^z_{\min} + n \frac{p-q}{2} - p + \frac{r + 11}{r-3} \opnorm{A - \E A} \leq \frac{r-3}{r-1} n \frac{p - q}{2},
\]
which is strictly more permissive than \eqref{eq:sbm_cond_pqknown}.
\end{proof}

\subsection{Asymptotic probabilistic analysis for the SBM}
\label{sec:proof_sbm_asymp}
Now, we show how \Cref{cor:sbm_determ} implies \Cref{thm:sbm_asymp}.
The proof structure is quite similar to that for \Cref{thm:z2_er_bern} in \Cref{sec:z2_bern}.
We will show that, with high probability, the condition \eqref{eq:sbm_cond} is satisfied (in which case the alternative condition \eqref{eq:sbm_cond_pqknown} also holds).

First, we show that the operator norm term in \eqref{eq:sbm_cond} is asymptotically negligible.
The condition \eqref{eq:sbm_asymp_cond} ensures $p \geq \frac{\log n}{n}$,
so we can apply \Cref{lem:X_conc} from \Cref{sec:z2_bern} with $v = p$ (note that $\E (A_{ij} - \E A_{ij})^2 \leq \E A_{ij}^2 \leq p$)
to obtain $\opnorm{A - \E A} \lesssim \sqrt{np}$ with probability $\to 1$ in $n$.
Furthermore, condition \eqref{eq:sbm_asymp_cond} implies the second inequality in
\[
	\frac{\sqrt{np}}{n(p-q)}
	\leq \frac{1}{\sqrt{n}} \frac{\sqrt{p} + \sqrt{q}}{p - q}
	= \frac{1}{\sqrt{n}(\sqrt{p} - \sqrt{q})} \\
	\leq \frac{1}{\sqrt{2 \log n}},
\]
so, with probability $\to 1$ in $n$,
\begin{equation}
	\label{eq:proof_sbm_opnorm}
	\frac{\opnorm{A - \E A}}{n(p - q)} \lesssim \frac{1}{\sqrt{\log n}}.
\end{equation}

Set, similarly to \eqref{eq:z2_proof_cs},
\begin{equation}
	\label{eq:sbm_proof_gamma_lb}
	\gamma = \parens*{\frac{1}{r-1} + \epsilon} (p - q) \quad \text{and} \quad \gamma' = \parens*{\frac{1}{r-1} + \frac{\epsilon}{2}} (p - q) = \gamma - \frac{\epsilon(p - q)}{2}.
\end{equation}
By \eqref{eq:proof_sbm_opnorm},
if we can show that $d^z_{\min} \geq \gamma' n$ with probability $\to 1$,
condition \eqref{eq:sbm_cond} will be satisfied with probability $\to 1$ as $n \to \infty$.

Note that $d^z_{\min}$ is the minimum of $n$ i.i.d.\ random variables of the form $V - W$,
where $V \sim \binomdist((n/2) - 1, p)$ and $W \sim \binomdist(n/2, q)$.
We focus on the lower tail of $V - W$, or, equivalently, the upper tail of $W - V$.
Note that, for any $t \in \R$,
\[
\E e^{t(W - V)} = (1 - q + q e^t)^{n/2}(1 - p + p e^{-t})^{\frac{n}{2} - 1},
\]
and, therefore,
\begin{align*}
	\log [\E e^{t(W - V)}]
	&= \frac{n}{2} \log(1 - q + q e^t) + \parens*{\frac{n}{2} - 1} \log (1 - p + p e^{-t}) \\
	&\leq \frac{n}{2} (-q + q e^t) + \parens*{\frac{n}{2} - 1} (- p + p e^{-t}) \\
	&= \frac{n}{2} ( - p + p e^{-t} - q + q e^t ) + p(1 - e^{-t}) \\
	&\leq \frac{n}{2} \brackets*{ p e^{-t} + q e^t - p - q } + 1.
\end{align*}
By a union bound and a Chernoff bound, for any $t \geq 0$,
\begin{align*}
	\P(d^z_{\min} \leq \gamma' n )
	&\leq n \P( W - V \geq -\gamma' n ) \\
	&\leq n \frac{\E e^{t(W - V)}}{e^{-t \gamma' n}} \\
	&\leq \exp \parens*{ \log n + \frac{n}{2} \brackets*{ p e^{-t} + q e^t - p - q} + \gamma' n t + 1} \\
	&= \exp \parens*{ \log n + n \brackets*{ \frac{p e^{-t} + q e^t}{2} - \frac{p + q}{2} + \gamma t} - \frac{\epsilon(p - q)}{2} n t + 1},
\end{align*}
where the last equality uses \eqref{eq:sbm_proof_gamma_lb}.
Noting that, for all $t \geq 0$, $t \leq \frac{e^t - e^{-t}}{2}$, we have
\begin{align*}
	\frac{p e^{-t} + q e^{t} }{2} + \gamma t
	&\leq \frac{(p - \gamma) e^{-t} + (q + \gamma) e^t}{2}.
\end{align*}
If $\gamma \leq \frac{p-q}{2}$ (which is ensured by the assumption $\epsilon \leq 1/12$), this last expression has a minimum value (over $t \geq 0$) of $\sqrt{(p - \gamma)(q + \gamma)}$
for $t = \frac{1}{2} \log \frac{p - \gamma}{q + \gamma}$.
Then, with this choice of $t$, we have
\begin{align*}
	\frac{p e^{-t} + q e^{t} }{2} - \frac{p + q}{2} + \gamma t
	&\leq \sqrt{(p - \gamma)(q + \gamma)} - \frac{p + q}{2} \\
	&= -\frac{1}{2} \parens*{ \sqrt{ p - \gamma } - \sqrt{q + \gamma} }^2.
\end{align*}
The condition \eqref{eq:sbm_asymp_cond} then ensures that
\[
	\log n + n \brackets*{ \frac{p e^{-t} + q e^t}{2} - \frac{p + q}{2} + \gamma t}
	\leq \log n - \frac{n}{2} \parens*{ \sqrt{ p - \gamma } - \sqrt{q + \gamma} }^2
	\leq 0,
\]
so we obtain
\[
	\P(d^z_{\min} \leq \gamma' n ) \leq \exp\parens*{ 1 - \frac{n\epsilon(p-q)}{4} \log \frac{p - \gamma}{q + \gamma} }.
\]
Denoting $\epsilon_r = \frac{1}{r-1} + \epsilon$ so that $\gamma = \epsilon_r(p-q)$,
we have
\begin{align*}
	\frac{1}{2} \log \frac{p - \gamma}{q + \gamma}
	&= \frac{1}{2} \log \frac{ \frac{p+q}{2} + (1-2\epsilon_r) \frac{p-q}{2}}{\frac{p+q}{2} - (1-2\epsilon_r) \frac{p-q}{2}} \\
	&\geq \frac{2}{p+q} (1 - 2 \epsilon_r) \frac{p-q}{2} \\
	&= (1 - 2 \epsilon_r) \frac{p-q}{p+q}
\end{align*}
by the inequality $\frac{1}{2} \log \frac{1 + x}{1-x} \geq x$.
Recalling that $r \geq 4$ and $\epsilon \leq \frac{1}{12}$,
we have $1 - 2 \epsilon_r \geq \frac{1}{6}$.
Furthermore, note that
\begin{align*}
	\frac{(p-q)^2}{p+q}
	&\geq \parens*{ \frac{p-q}{\sqrt{p} + \sqrt{q}} }^2 \\
	&= (\sqrt{p} - \sqrt{q})^2 \\
	&\geq (\sqrt{p-\gamma} - \sqrt{q+\gamma})^2 \\
	&\geq \frac{2 \log n}{n}
\end{align*}
by the condition \eqref{eq:sbm_asymp_cond}, so, putting these inequalities together, we obtain
\begin{align*}
	\frac{n\epsilon(p-q)}{4} \log \frac{p - \gamma}{q + \gamma}
	&\geq \frac{n\epsilon(p-q)}{2}(1 - 2 \epsilon_r) \frac{p-q}{p+q} \\
	&\geq \epsilon(1 - 2 \epsilon_r) \log n \\
	&\geq \frac{\epsilon}{6} \log n,
\end{align*}
so
\[
	\P(d^z_{\min} \leq \gamma' n ) \leq \exp\parens*{ 1 - \frac{\epsilon}{6} \log n },
\]
which goes to zero as $n \to \infty$.
Thus we obtain the result \Cref{thm:sbm_asymp}.
	
\section{Robustness to monotone adversaries}
\label{sec:monotone}
A useful property for synchronization and clustering algorithms is robustness to a ``monotone adversary.''
For example, for the synchronization problem of \Cref{sec:intro_er},
what happens if somebody is allowed to add additional edges to the measurement graph with correct measurements or to delete or correct some of the erroneous measurements?
For the graph clustering problem of \Cref{sec:intro_sbm},
what happens if somebody adds edges between vertices in the same ground-truth cluster and/or removes edges between vertices in different clusters?
With the random base models we have presented, adding the possibility of additional such ``helpful'' modifications gives a \emph{semi-random} model \cite{Feige2001}.

Intuitively, such an ``adversary'' should only be helping us and making the problem easier.
However, algorithms that depend on certain regularity properties (e.g., matrix eigenvalues/vectors) of the underlying random model may be severely disrupted.
It is in fact possible that the underlying estimation problem becomes fundamentally \emph{harder}. See, for example, \citet{Moitra2016} for further discussion and references.
Nevertheless, certain problem settings and algorithms do exhibit robustness.

In our setting, we characterize this robustness as follows.
Suppose, given a cost matrix $C$, the problem \eqref{eq:obj_orig} has optimal solution $z$.
Now suppose we perturb $C$ by a matrix $\Delta^+$ such that, for all $i,j$,
\begin{equation}
	\label{eq:monotone_deltaplus}
	\Delta^+_{ij} z_i z_j \geq 0.
\end{equation}
Clearly, \eqref{eq:obj_orig} with cost matrix $C' = C + \Delta^+$ will still have $z$ as an optimum (and, if the optimum is unique for the original problem, it also is for the perturbed problem).
An algorithm is robust to monotone adversaries if, given that it successfully recovers the optimum $z$ from data $C$, it \emph{also} recovers $z$ with data $C' = C + \Delta^+$ for \emph{any} $\Delta^+$ satisfying \eqref{eq:monotone_deltaplus}.

One can easily verify (see, e.g., \cite{Feige2001,Hajek2016}) that the SDP relaxation \eqref{eq:obj_sdp} of \eqref{eq:obj_orig} has such robustness.
What about the nonconvex partial relaxation \eqref{eq:obj_bm}?
\citet{Ling2019} showed, by a connection to Kuramoto oscillator networks, that in the case $r = 2$ the nonconvex approach is \emph{not} robust to monotone adversaries.

However, when $r \geq 3$, \citet{Ling2023b} showed that the nonconvex approach is robust in the sense that the benign landscape results in that paper are not harmed by monotone adversaries.
The same is true for our results:
\begin{theorem}
	\label{thm:monotone}
	The results \Cref{thm:z2_determ,thm:z2_er_bern} and, in the case where $C$ is set by \eqref{eq:sbm_C_pqknown}, \Cref{cor:sbm_determ,thm:sbm_asymp} are robust to monotone adversaries in the following sense:
	for each theorem, under the same conditions, if we further perturb $C$ by replacing it with $C' = C + \Delta^+$,
	where $\Delta^+$ and $z$ satisfy \eqref{eq:monotone_deltaplus},
	the same conclusion holds.
\end{theorem}
To see this, note that in the context of \Cref{thm:z2_determ}, we can replace the graph adjacency matrix $A$ with the matrix $A'$ defined by $A'_{ij} = A_{ij} + \Delta^+_{ij} z_i z_j$.
Because each entry is only increased, the algebraic connectivity $\lambda_2$ is only increased.
This shows robustness for \Cref{thm:z2_determ}; the same result immediately propagates to that theorem's direct corollaries.

For the ordinary graph clustering problem with $C$ chosen as in \eqref{eq:sbm_C},
the situation is less clear,
as elementwise perturbations affect the entire cost matrix $C$ via the centering operation.
There are certainly ways around this issue for the full SDP relaxation (see \cite{Hajek2016,Bandeira2018}), but an extension to our nonconvex setting is not trivial, and we do not pursue it in this paper.

\Cref{thm:monotone} does \emph{not} say that for \emph{every} cost matrix $C$ such that \eqref{eq:obj_bm} happens to have a benign landscape, a monotone perturbation will preserve that landscape.
That would be a stronger result and would require quite different proof techniques (as ours fundamentally depend on the problem structure).
However, the specific analysis that leads to the result \Cref{thm:z2_determ} and its corollaries is only helped by such perturbations. 

\section*{Acknowledgements}
This work was supported by the Swiss State Secretariat for Education, Research and Innovation (SERI) under contract MB22.00027.

\bibliography{refs}

\end{document}